\documentclass{article}

\usepackage{arxiv}

\usepackage[utf8]{inputenc} 
\usepackage[T1]{fontenc}    
\usepackage{hyperref}       
\usepackage{url}            
\usepackage{booktabs}       
\usepackage{amsfonts}       
\usepackage{nicefrac}       
\usepackage{microtype}      
\usepackage{lipsum}		
\usepackage{graphicx}
\usepackage{natbib}
\usepackage{doi}
\usepackage{amsmath}

\newcommand\figref[1]{\textbf{Figure~\ref{fig:#1}}}
\newcommand\tabref[1]{\textbf{Table~\ref{tab:#1}}}

\title{English version of "Screen-line counter location problem with O/D cut selection approach"}

\author{ \href{https://orcid.org/0000-0003-4341-2155}{\includegraphics[scale=0.06]{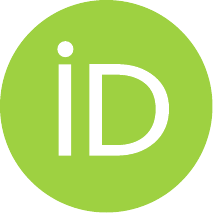}\hspace{1mm}Satoshi Sugiura}\thanks{This secondary publication is an English translation of doi:(https://doi.org/10.2208/jscejj.23-00244)} \\
	Graduate School of Engineering\\
	Hokkaido University\\
	Kita 13, Nishi 8, Sapporo, Hokkaido \\
	\texttt{sugiura@eng.hokudai.ac.jp} \\
}

\renewcommand{\headeright}{Secondary  publication}
\renewcommand{\shorttitle}{SCLP with O/D cut selection approach}

\hypersetup{
pdftitle={Screen-line counter location problem with O/D cut selection approach},
pdfsubject={q-bio.NC, q-bio.QM},
pdfauthor={Satoshi Sugiura},
pdfkeywords={traffic sensor location problem, screen-line counter location problem, graph cut, maximum weighted closure problem, mixed integer linear programming problem},
}

\begin{document}
\maketitle

\begin{abstract}
  This paper provides an efficient solution approach to the screen-line counter location problem (SCLP), which is a counter location problem with the constraint that the traffic between OD pairs must be observed at least once.
  This paper formulates the SCLP using a graph cut approach, which consists of an enumeration of cuts and a cut selection problem.
  These problems can be reduced to a concise formulation that extends the maximum weight closure problem for two problems: finding the minimum number of links that observe all OD pairs and finding the maximum number of observed OD pairs with a budget-constrained number of links.
  Insights into the characteristics of cuts give superior upper bounds on the problem of finding the minimum number of links that observe all OD pairs.
  The proposed method is evaluated on the Sioux-Falls network. It shows that it is possible to derive a solution equivalent to the optimal solution found in previous studies in a very short computation time.
\end{abstract}

\keywords{traffic sensor location problem \and screen-line counter location problem \and graph cut \and maximum weighted closure problem \and mixed integer linear programming problem}

\section{Introduction}
Traffic data is essential for formulating and planning traffic policies.
It is primarily collected from traffic sensors installed on the road network. 
These sensors are crucial in understanding traffic situations, including traffic volume and route choice conditions. 
Due to the importance of the acquired data, there is significant interest in determining the optimal locations for installing traffic sensors. 
To observe traffic conditions between regions, the most empirical method of selecting locations that efficiently capture traffic conditions on the network by passing through mountainous areas, overpasses, and river bridges has been used.

Traffic sensor location problem (TSLP) deals with the optimal location of traffic sensors according to each purpose of traffic state observation.
TSLP has evolved as a methodology according to the purpose of observation and the type of sensor, starting with \citet{Yang1998}.
\citet{Owais2022_review} provided a literature review related to TSLP over the past 30 years and classified the topics into the following six categories according to the purpose of observation.
\begin{enumerate}
  \item OD trip table estimation/updating \citep{Yang1998, AChen_2007}, 
  \item flow observability \citep{Castillo_2008, Bianco_2014}, 
  \item link flow inference \citep{Hu_2009, Ng_2013},
  \item path reconstruction \citep{Gentili_2005, Minguez_2010},
  \item screen line/traffic surveillance \citep{Yang_2001,Yang2006},
  \item travel time estimation \citep{Viti_2008, Mirchandani_2009},
\end{enumerate}

This paper focuses on the Screen-line counter location problem (SCLP)\citep{Yang2006}, which is related to the fifth category.

\citet{Yang1991} emphasized the importance of considering the links to be observed as input for OD estimation and showed the basic principle of sensor location called OD Covering Rule.
This rule states that for the Maximum possible relative error of the estimated OD traffic volume to take a finite value, sensors must be placed so that each OD pair is observed at least once.
Subsequently, \citet{Yang1998} showed four basic sensor location rules: OD Covering Rule, maximal flow fraction rule, maximal flow-intercepting rule, and link independence rule.
Many researchers have studied TSLP based on these rules.
SCLP is a TSLP that optimizes sensor location with the OD Covering Rule as a constraint.
\citet{Yang2006} have shown two problems of SCLP: one is to find the minimum number of links to observe all OD pairs, and the other is to find the maximum number of observable OD pairs with a budget-constrained number of locating sensors.
SCLP was initially referred to as a screen-line counter location problem. However, in practice, it has been approached as the problem of identifying links that meet the two criteria defined.
\citet{Owais2019} applied SCLP to determine the sensor location strategy to detect abnormal vehicles emitting exhaust gases.
SCLP originates from OD estimation and has various applications for sensor location to observe trips between OD pairs.

The previous method for SCLP required repeated path enumeration for column generation, leading to increased computational load.
The reason is that as the number of paths increased, it became difficult to search for the exact optimal solution by combining them.
Therefore, heuristic methods have been used to reach the solution \citep{Yang2006, Owais2019, Owais2022_IEEE}.
However, these heuristic methods limit the number of paths to be enumerated to reduce the computational load, leaving the problem of remaining paths that do not have a counter between $(s,t)$ \citep{Owais2022_IEEE}.
It is challenging to construct a method that guarantees the satisfaction of the OD covering rule and has a low computational load.

This study aims to achieve the OD coverage rule by utilizing the characteristics of graph cuts and providing a method to solve SCLP with a lower computational load than previous studies.
The $(s,t)$-cut always separates the graph for any $(s,t)$, which means that at least one of the links on the $(s,t)$-cut must pass through any path between $(s,t)$.
This means that the set of sensors captures all paths between $(s,t)$ if these are located on all links included in the cut.
This paper provides a new method for SCLP with enumerating cuts and identifying the optimal counter location links by selecting the enumerated cuts.
Since $(s,t)$-cut can be enumerated more efficiently in advance than paths between $(s,t)$, this study has the advantage of not requiring sequential enumeration of paths or iterative calculations, as in conventional SCLP such as \citet{Yang2006}.
The proposed method is characterized by the fact that the optimal solution can be obtained by executing the cut enumeration algorithm and the optimization problem only once. 

The main contributions of this paper are as follows.
\begin{itemize}
  \item Describing SCLP using a graph cut approach.
  \item Providing a solution approach consisting of cut enumeration and a cut selection problem (CSP). CSP is reduced to a concise formulation that extends the maximum weighted closure problem (MWCP).
  \item Showing that CSP corresponds to two screen-line counter location problems proposed by \citet{Yang2006}: finding the minimum number of links to observe all OD pairs and finding the maximum number of observable OD pairs with a budget-constrained number of links.
  \item Providing a superior upper bound on the problem of finding the minimum number of links to observe all OD pairs and a method to limit the cuts to be enumerated based on this upper bound.
  \item The proposed method is verified for performance on the Sioux-Falls network and shows that a solution equivalent to the optimal solution found in previous studies can be obtained in a very short computation time.
\end{itemize}

\section{SCLP formulation based on cut property}\label{sec:formulation}

\subsection{Network description}

Provided that $G(V,E)$ is a directed road network, where $V$ is the node set and $E(=\{\forall(i,j)|i,j\in V,(i,j)\neq (j,i)\})$ is the set of directed links.
$Q(\subseteq V)$ is the set of all centroids (origin and destination points) on the network, and $w=(s,t)\mid s,t\in Q$ is an OD pair, where $W$ is the set of OD pairs.
The OD pair is denoted as $(s,t)$ throughout this paper.
There are several definitions of cut, but this paper employs the following definition of cut. 

\textbf{Definition: $(s,t)$-cut} A $(s,t)$-cut is a set of directed links that, when removed from the directed graph $G(V,E)$, eliminates the existence of a directed path from $s$ to $t$.
$C^{st}$ is the set of all $(s,t)$-cuts, and its components $c_m^{st}$ are constructed by links according to the definition ($c_m^{st} \subseteq E$).
$M(st)$ is the total number of $(s,t)$-cuts.
The number of links included in the cut $c_m^{st}$ is denoted as $|c_m^{st} |$.
$P_w$ is the set of all simple directed paths between $(s,t)=w$.
$\Xi_{wk}^e$ is the link-path incidence matrix that takes 1 when link $e$ is included in the $k$-th path of the path set $P_w$ and 0 otherwise.

\subsection{Advantages of cut-based formulation}
The basic formulations of SCLP that have been constructed so far \citet{Yang2006} are shown below.
Let $x_e$ be a binary variable that takes 1 when a counter is placed on link $e$ and 0 otherwise.
First, the formulation to find the minimum number of links to observe all OD pairs is shown below.
\begin{equation}
  \label{eq:sclp_1_obj}
  \min_{\mathbf{x}}{\sum_{e\in E}{x_e}}
\end{equation}
\begin{equation}
  \label{eq:sclp_1_cns}
  \sum_{e\in E}{\Xi_{wk}^e x_e}\geq 1 \quad \forall k\in P_w, w \in W
\end{equation}
\begin{equation}
  \label{eq:sclp_1_defval}
  x_e = \{0,1\} \quad \forall e \in E
\end{equation}
Eq.\eqref{eq:sclp_1_obj} is the objective function that minimizes the number of counter-placed links.
Eq.\eqref{eq:sclp_1_cns} is the constraint condition that at least one link with a counter is placed in all paths in $w$.
This constraint ensures that all $w\in W$ and all paths are observed at least once.

Below is the basic formulation for finding the maximum number of OD pairs with a budget-constrained number of links.
Let $\alpha_w$ be 1 when $(s,t)=w$ is not observable, and 0 otherwise.
\begin{equation}
  \label{eq:sclp_2_obj}
  \max_{\mathbf{x}}{\sum_{w\in W}{\alpha_w}}
\end{equation}
\begin{equation}
  \label{eq:sclp_2_cns}
  \sum_{e\in E}{\Xi_{wk}^e x_e}\geq 1 \quad \forall k\in P_w, w \in W
\end{equation}
\begin{equation}
  \label{eq:sclp_2_defval}
  x_e = \{0,1\} \quad \forall e \in E
\end{equation}
Eq.\eqref{eq:sclp_2_obj} is the objective function that maximizes the number of observable OD pairs with a budget-constrained number of links.
Eq.\eqref{eq:sclp_2_cns} is the constraint condition that at least one link with a counter is placed in all paths in $w$.
This constraint ensures that all $w\in W$ and all paths are observed at least once.

These formulations require enumerating all simple directed paths between $(s,t)$ for $w\in W$.
However, the number of paths increases exponentially with the number of nodes and links in the network, making it difficult to find the exact optimal solution by combining them.
The previous methods for SCLP have been based on the column generation method.
The solution methods consist of two steps: finding the optimal sensor location for the enumerated paths and enumerating the remaining paths that do not have a counter between $(s,t)$.
Specifically, $\dot{\Xi}_{wk}^e$ is a link-route incidence matrix based on the condition in the counter installation link decision problem, the subset $\dot{P}_w$ of the set of all routes $P_w$, and solve the limited master problem with this as input.
The limited master problem is solved by searching for new routes based on the dual cost obtained by the subproblem.
These two steps are repeated until no new route with a small dual cost is found.

\begin{figure}[t]
  \centering
  \includegraphics[keepaspectratio,scale=0.6]
  {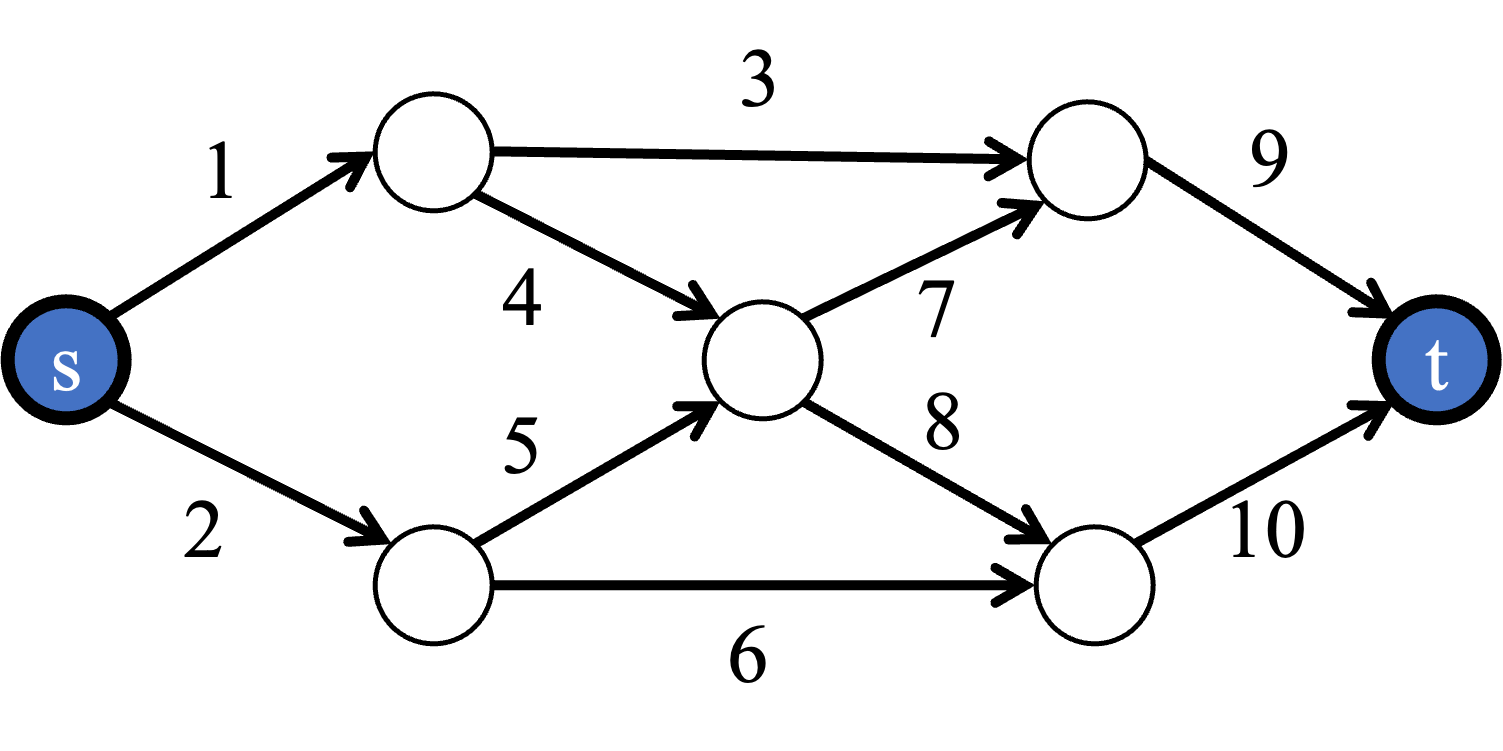}
  \caption{The sample network.}
  \label{fig:ig_network}
\end{figure}

Here, I explain the difficulty of the conventional column generation method with \figref{ig_network}.
For simplicity, the network is a simple network with one origin and destination nodes $(s,t)$.
The numbers in the figure represent the link ID.

The subproblem in the column generation method provides the link costs. 
If the sensor is located on the link, the cost is 1; otherwise, it is 0. In this situation, the shortest path is searched.
Initially, in this example, no counters are placed on any links, so arbitrary paths are searched. Let's assume that the path $\{1,4,7,9\}$ is obtained.
Due to the restricted master problem, I select a minimal set of links for each listed route, ensuring that each route contains at least one link.
For example, the sensor can be located on the link 4.
This process is repeated until no new routes with 0 cost are found.
Assume that the next route search results in a path of $\{2,6,10\}$ that does not include link $4$ and that the counter is located at $\{4,10\}$.
Then the third iteration searches for a path that does not contain $\{4,10\}$.
These iterations are repeated until no new routes with zero cost are found.
In this iterative process, if the set of links $\{1,2\}$ or $\{9,10\}$ are selected, the optimal solution is obtained and the iteration is terminated.
Note that routes are searched based on cost, with the only requirement being to avoid links with installed counters.
However, as the costing clearly shows, routes are searched, with the only requirement being to avoid links where counters are placed.
For example, a subgraph region consisting of white nodes will have many alternative paths, and new paths will be discovered one after another, even if a small number of counters are placed on these links.
It is understood that for complex networks with many different paths, the number of iterations can be huge, as the optimal solution is reached by chance or, at worst, it is iterated until all paths are enumerated.

An approach based on cut is proposed to solve this issue mentioned.
Any $(s,t)$-cut necessarily separates the node pair $(s,t)$, which means that at least one of the links on the $(s,t)$-cut must pass through any path between $(s,t)$.
This means that the set of sensors captures all paths between $(s,t)$ if these are located on all links included in the cut $c^{st}\in C^{st}$.
$\{1,2\}$ and $\{9,10\}$ are both the $(s,t)$-cut and there are no path between $(s,t)$ without these links.
This means that the proposed method in this paper can employ the $(s,t)$-cut to find the optimal solution satisfying OD Covering Rule.

Therefore, the proposed method is to enumerate all $(s,t)$-cuts and select the optimal counter location links by selecting the enumerated cuts.

\subsection{The method to enumerate the cut}\label{sec:cut_enu}
\citet{provan_1996} provided a method to enumerate all $(s,t)$-cuts in the directed graph.
They showed that the time-per-cut complexity for the network size was linear and could be calculated efficiently.
There are other algorithms to enumerate cuts than the method proposed by \citet{provan_1996} such as ZDD(Zero-suppressed Binary Decision Diagram)\citep{kawahara_2012}
This paper employs the algorithm proposed by \citet{provan_1996} to enumerate all $(s,t)$-cuts.

\subsection{Formulation as the non-linear integer programming problem}
There exist multiple $(s,t)$-cuts to separate $(s,t)\in W$.
When the sensors are located in all links that include one of these cuts, all paths between $(s,t)$ are observed.
Therefore, the cut chosen from the enumerated $(s,t)$-cuts should be selected to optimize the objective function.
Based on this concept, a simple formulation involving the multiplication of unknown variables is initially considered.
The formulation of the problem given in \citet{Yang2006}, the problem of giving the minimum number of counter location links under the condition of observing between all $(s,t)$ is called (NCSP1), and the problem of finding the maximum number of $(s,t)$ observations with a constraint on the number of counter placement links is called (NCSP2).
I define the following variables. $x_{ij}$ is a binary variable that takes 1 when the sensor is located on the link $(i,j)$, and 0 otherwise.
$\gamma_{m}^{st}$ is a binary variable that takes 1 when the cut $c_m^{st}$ is selected, and 0 otherwise.
$\delta(\cdot)$ is the step function that takes 0 when the input equals zero and 1 when the input is larger than 0.
$\sigma_{ij}(c_m^{st})$ takes 1 when the link $(i,j)$ is included in the cut $c_m^{st}$, and 0 otherwise.

NCSP1 is formulated as follows.
\begin{equation}
  \label{eq:NCSP1_obj}
  \min_{\mathbf{x}}{\sum_{(i,j)\in E}{x _{ij}}} 
\end{equation}
\textit{subject to}\\
\begin{equation}
  \label{eq:NCSP1_cut_number}
  \sum_{s,t,m \mid c_m^{st} \in C^{st}}{\prod_{(i,j)\in c_m^{st}}{\sigma_{ij}{(c_m^{st})x_{ij}}}}\geq 1 \ \forall st|(s,t)\in W
\end{equation}
\begin{equation}
  \label{eq:lambda_def}
  x_{ij}= \{0,1 \} \ \forall (i,j) \in E
\end{equation}
Eq.\eqref{eq:NCSP1_obj} is the objective function that minimizes the number of counter-locating links.
Eq.\eqref{eq:NCSP1_cut_number} is the constraint condition that at least one of the cut $c_m^{st}$ is selected and located sensor on all links included in the cut.

NCSP2 is formulated as follows. $K$ is the upper limit of the number of counter-locating links.
\begin{equation}
  \label{eq:NCSP2_obj}
  \max_{\mathbf{x}}{\sum_{st|(s,t) \in W}{\delta \left ( \sum_{s,t,m \mid c_m^{st} \in C^{st}}{\prod_{(i,j)\in c_m^{st}}{\sigma_{ij}{(c_m^{st})x_{ij}}}}\right )}}
\end{equation}
\textit{subject to}\\
\begin{equation}
  \label{eq:NCSP2_link_number}
  \sum_{(i,j)\in E}{x_{ij}}\leq K
\end{equation}
\centerline{Eq.\eqref{eq:lambda_def}}\\

Eq.\eqref{eq:NCSP2_obj} is the objective function to maximize the number of separated $(s,t)$ pairs. 
Eq.\eqref{eq:NCSP2_link_number} is the budget constraint of the number of locating sensors.

The above formulations are non-linear integer programming problems involving the multiplication of step functions and variables, and both are difficult to solve.
If a solution algorithm for this formulation is provided in the future, it can be solved as is.
In the following sections, two CSPs are formulated in forms that are easier to solve than NCSP1 and NCSP2.

\subsection{Formulation as an extension of Maximum Weighted Closure Problem}
\subsubsection{MWCP overview}
\begin{figure}[t]
  \centering
  \includegraphics[keepaspectratio,scale=0.5]
  {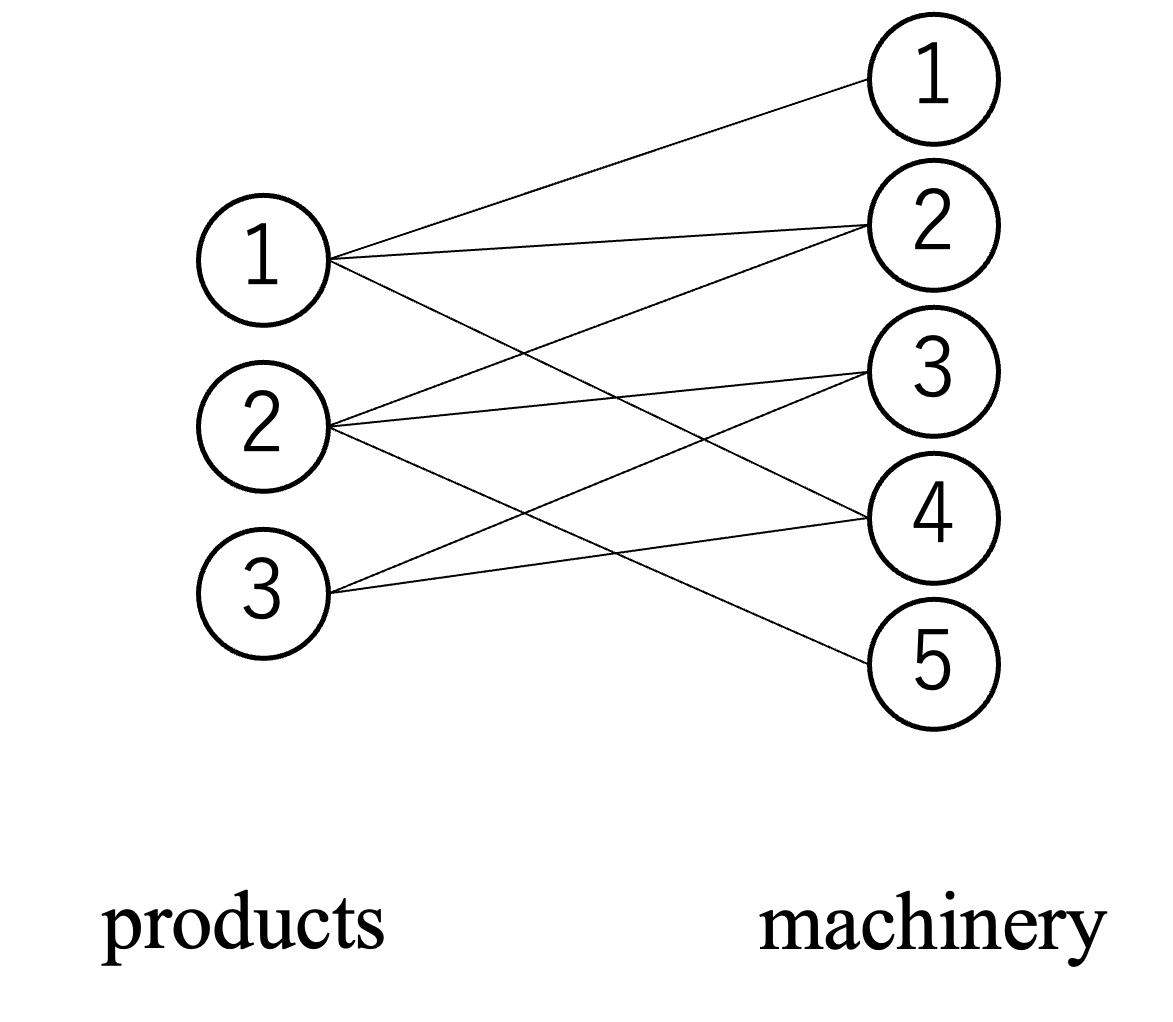}
  \caption{An example of MWCP}
  \label{fig:eg_mwcp}
\end{figure}

The Maximum Weighted Closure Problem (MWCP) involves selecting the optimal subset for two different sets.
One example is shown in \figref{eg_mwcp}.
Assuming I have a set of products under development and a set of machine tools required to produce those products.
The set of products is represented by $N^L=\{1,2,3\}$, and the benefit to make each product is set as $u_l (l\in N^L)$.
The set of tools is represented by $N^R=\{1 ,2,3,4,5\}$, and the cost of using each tool is set as $t_r (r\in N^R)$.
Each product requires a specific set of tools; if any of the needed tools are lacking, the product cannot be produced.
This relationship of products and tools is represented by the bipartite graph.
$D(N,A)$ is the bipartite graph where $N\left(=N^L\bigcup N^R\right)$ is the set of nodes of this bipartite graph and $A$ is the set of arcs.

Under this condition, the administrator wants to solve the sets of products and tools that maximize the sum of the benefits (not negative) and the costs (not positive).
Note that in this example, the amount of product manufactured, the cost of materials required for the product, etc. are discarded.
Assuming that the decision variables indicating the choice of product $l$ are $y_l=\{0,1\}$ and the decision variables indicating the choice of machine tool $r$ are $x_r=\{0,1\}$, the MWCP can be formulated as follows.
\begin{equation}
  \label{eq:obj_orig}
  \max_{\mathbf{y,x}}{\sum_{l\in N^L}{u_l y_l}-\sum_{j\in N^R}{t_r x_r}} 
\end{equation}
\textit{subject to}\\
\begin{equation}
  \label{eq:cut_link_cns}
  y_l \leq x_r \quad \forall (l,r)\in A 
\end{equation}
\begin{align}
  y_l = \{0,1\} \quad \forall l\in N^L\\
  x_r = \{0,1\} \quad \forall r\in N^R
\end{align}
Eq.\eqref{eq:obj_orig} is the objective function of MWCP, which maximizes the sum of the benefits of product manufacturing and the total cost of machine tools.
Eq.\eqref{eq:cut_link_cns} represents the relationship between product manufacturing and machine tool procurement.

Constrain that each corresponding machine tool is procured ($x_r=1$) to manufacture the product ($y_l=1$).
Since the left-hand side of the constraint is a fully unimodular matrix, it is known that a linear relaxation to $0\leq y_l,x_r\leq 1$ yields a solution in the form {0,1}.
The $\mathbf{x,y}$ obtained by solving the optimization problem represents a subset of the objective function from the different sets of products and machine tools.
In other words, in this example, the output is the set of products to be manufactured and the set of machine tools needed to manufacture them.

\subsubsection{Generalized formulation as an extension of MWCP}

This paper considers a set of cuts and their constituent links as a set and formulates them as an extension of MWCP by describing their relationship in a bipartite graph.
To prepare the formulation, the enumerated cuts $c_m^{st}\in C^{st} \forall st|(s,t)\in W$ and their corresponding links as matrices are represented.
That is to say, define a matrix $\Psi$ with size $\left (\sum_{st|(s,t)\in W} M(st) ,|E|\right )$ and let $\psi_{lr}$ be its component. Let $l$ be the row, and $r$ be the column identifier.
Every enumerated cut $c_m^{st} \in C^{st} \forall st|(s,t)\in W$ is given a uniquely identifying number, as are the links $e^\in E$.
The cut and link identifiers $l,r$ are output by the functions $\zeta,\eta$, which are expressed as $\zeta(c_m^{st} )=l,\eta(i,j)=r$, respectively. The $\psi_{lr}$ is defined by Eq.\eqref{eq:psi}.

\begin{align}
  \psi_{lr}=
  \begin{cases}
    1 & if \zeta (c_m^{st})=r, \eta (i,j)=l |\\
    &  (i,j)\in c_m^{st},c_m^{st}\in C^{st}, st|(s,t)\in W\\
    0 & otherwise
  \end{cases}
  \label{eq:psi}
\end{align}

From the generated matrix $\Psi$, a bipartite graph showing the relationship between the cut $c_m^{st}$ and the corresponding link $(i,j)\in E$ can be constructed.

$N^L$ and $N^R$ are $|N^L|=\sum_{st|(s,t)\in W}{M(st)}$ and $|N^R|=|E|$, consistent with the row and column sizes of $\Psi$, respectively.
The link $(l,r)\in A|l\in N^L, r\in N^R$ connects only node pairs with $\psi_{lr}=1$.
This bipartite graph shows $N^l$-enumerated cuts as nodes, and the node corresponding to each cut is connected by a link to the node $r\in N^R$ corresponding to the link on the road network contained in the cut.

\begin{figure}[tb]
  \centering
  \includegraphics[keepaspectratio,scale=0.5]
  {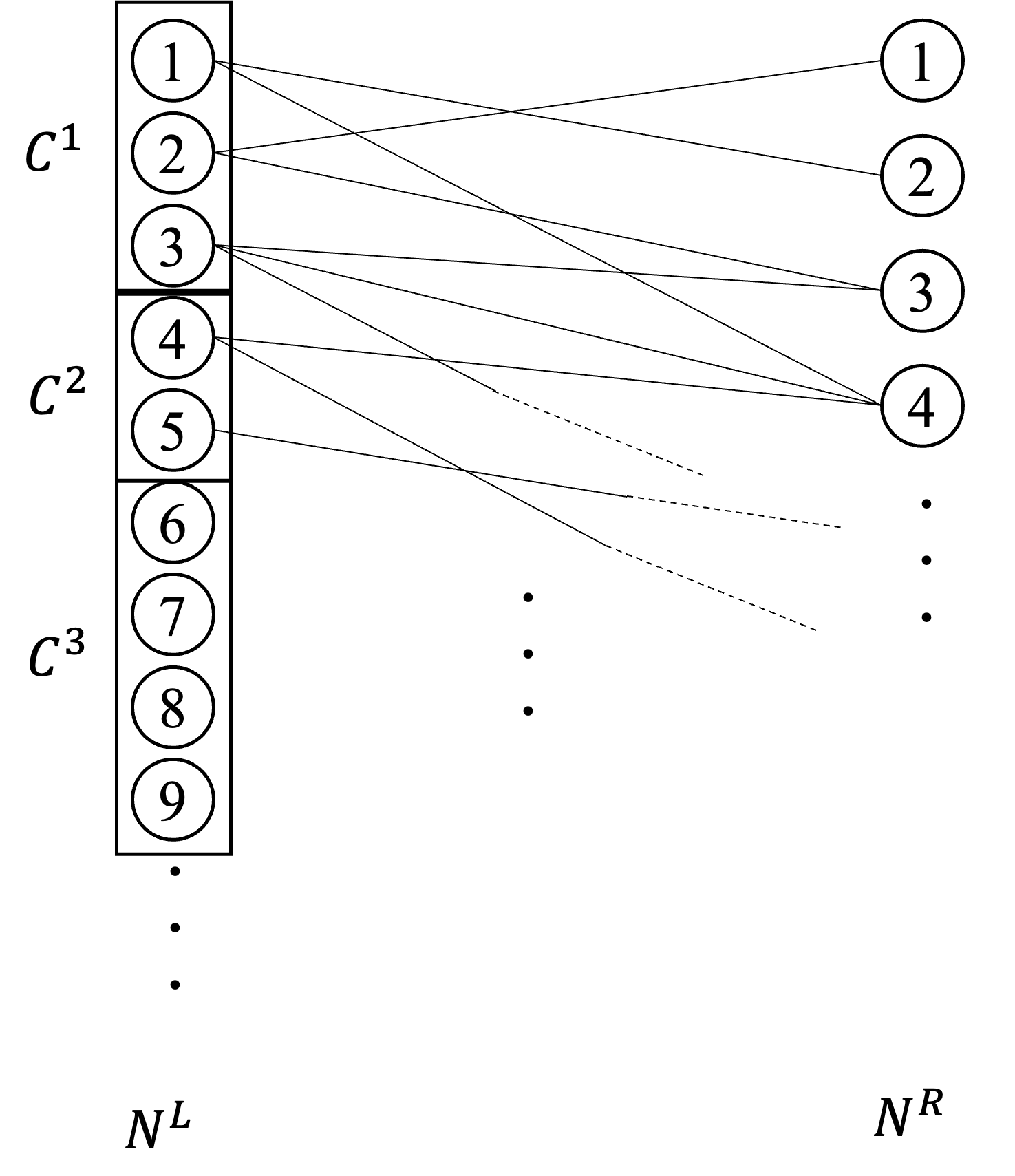}
  \caption{Example of a bipartite graph constructed by a matrix $\Psi$}
  \label{fig:bipartite}
\end{figure}

The policy of applying MWCP to the problem addressed in this paper is illustrated in the example of \textbf{Figure \ref{fig:bipartite}}.
Each node on the left side of the bipartite graph corresponds to a cut $c_m^{st}$, and it is possible to identify which $(s,t)$-cut it is.
In the figure, for example, the cut contained in $C^1 \left ((s,t)=w=1\right )$ is represented by node $\{1,2,3\}$.
If the cut indicated by node $l$ is selected, a profit $u_l$ is obtained, but it costs $t_r$ to place a counter on the link indicated by node $r$.

The decision to employ a cut $l$ that observes between $(s,t)$ is represented by $y_l=\{0,1\}$, and the location of a counter on link $r$ is represented by $x_r=\{0,1\}$.
Counters must be installed on all links that make up the cut when employing a cut that observes between $(s,t)$.
Therefore, the relationship between the employment of a cut and the decision to place counters on links can be expressed in the same way as in the MWCP example by the Eq.\eqref{eq:cut_link_cns}.
Although both examples show the relationship between the employment of cuts and the placement of counters on links, Eq.\eqref{eq:cut_link_cns} describes the relationship as links in a bipartite graph, and the separation into constraints for each link allows for simplified expressions that avoid multiplication by unknown variables.
Under the constraints of Eq.\eqref{eq:cut_link_cns}, the objective function of MWCP can be expressed as in the MWCP example with Eq.\eqref{eq:obj_orig} that maximizes the difference between total benefits and total costs.

The objective function Eq.\eqref{eq:obj_orig}, constraint condition Eq.\eqref{eq:cut_link_cns}, $0\leq y_l,x_r\leq 1$ is the formulation of the cut selection problem as a basic MWCP. The optimal solution identifies the cut to be employed and the set of links necessary for its observation The optimal solution is determined by $\mathbf{y,x}$.

In the proposed method, multiple cuts to observe between the same $(s,t)$ are enumerated (i.e., $C^{st}=\{c_1^{st},c_2^{st},. .c_m^{st}. .c_M^{st}\}$).
Since it is inefficient in SCLP for the same $(s,t)$ to be observed by multiple cuts, and it is also unreasonable to duplicate gains, I constrain the following Eq.\eqref{eq:y_cns} so that each $(s,t)$ is employed by one or less cuts to be observed Constrain $(s,t)$ to adopt no more than one cut to observe between each of them.
This means that in \textbf{Figure \ref{fig:bipartite}}, for example, nodes $1,2,3$ are selected from $C^1$, from $C^2$, nodes $4,5$ are selected one or less from the corresponding cuts.

\begin{align}
  \label{eq:y_cns}
  \sum_{l=\zeta(c_m^{st}|c_m^{st}\in C^{st})}{y_l}\leq 1 \quad \forall st|(s,t)\in W 
\end{align}

Note that the Heaviside function used in NCSP2 is not necessary for the objective function in the subsequent formulations since this constraint expresses that if one $(s,t)$-cut is employed, the corresponding $(s,t)$-pair of traffic is observed.

The above Eqs.\eqref{eq:obj_orig}, \eqref{eq:cut_link_cns}, and {\eqref{eq:y_cns}} are the general formulation (CSP0) of the cut choice problem as an extension of the MWCP.
However, since CSP0 is no longer a fully unimodular matrix due to the Eq.\eqref{eq:y_cns}, $y_l$ requires an integer condition, making linear relaxation impossible.
In contrast, $x_r$ can be linearly relaxed to an integer as long as $\mathbf{y}$ is an integer.

Although CSP0 assumes that the profit $u_l$ from observing traffic between $(s,t)$ and the cost of counter installation $t_r$ are linearly combinable, in practice, it would be difficult to identify the parameters that combine them.
Considering practical conditions, formulation consistent with the two problems presented in \citet{Yang2006} is provided.

\subsubsection{Allocation problem for the minimum number of links to be observed between all $(s,t)$}
Let CSP1 be the problem of achieving counter allocation with minimum links under the condition that all $(s,t)\in W$ are observed.
CSP1 can be expressed by replacing the objective function as Eq.\eqref{eq:obj_orig} and constraint as Eq.\eqref{eq:y_cns} in CSP0 with the following Eqs.\eqref{eq:obj_csp1} and \eqref{eq:y_cns_csp1} respectively.

\begin{equation}
  \label{eq:obj_csp1}
  \min_{\mathbf{y,x}}{\sum_{j\in N^R}{t_r x_r}} 
\end{equation}
\begin{align}
  \label{eq:y_cns_csp1}
  \sum_{l=\zeta(c_m^{st}|c_m^{st}\in C^{st})}{y_l}= 1 \quad \forall st|(s,t)\in W 
\end{align}

Eq.\eqref{eq:y_cns_csp1} constrains the employment of only one observable cut between every $(s,t)\in W$.
Note that CSP1 can be linearly relaxed compared to CSP0 since the constraints related to $\mathbf{x}$ are identical.
Eq.\eqref{eq:y_cns_csp1} can be defined as a Special Order Set type1 (SOS1) constraint, which allows us to solve Brunch and Bound in binary integer programming problems more efficiently.

The problem can be solved as the number of links to install counters if $t_r$ are all the same value, and the total installation cost can also be considered by setting the installation cost according to the individual link characteristics.
CSP1 can be formulated as a mixed integer linear programming problem, which can be solved efficiently by using a standard solver, especially one that supports the SOS1 constraint.

\subsubsection{Location problem to maximize the number of observed $(s,t)$ pairs under a budget for the number of links}
Let CSP2 be the problem of observing as many $(s,t)$ pairs as possible under a budget for the number of links.
CSP2 can be expressed by replacing Eq.\eqref{eq:obj_orig} of CSP0 with the following Eq.\eqref{eq:obj_csp2} and adding the Eq.\eqref{eq:x_budget}.
\begin{equation}
  \label{eq:obj_csp2}
  \max_{\mathbf{y,x}}{\sum_{l\in N^L}{u_l y_l}} 
\end{equation}
\begin{equation}
  \label{eq:x_budget}
  \sum_{r\in N^R}{x_r}\leq K
\end{equation}
$K$ is the same definition as in NCSP, i.e., the budget for the number of links where the counter is installed.
However, CSP2 requires that $x_r=\{0,1\}$, since the solution is not promised to be an integer when linearly relaxed by the Eq.\eqref{eq:x_budget}.

If $u_l$ is identical for both cut $l$ and the corresponding $(s,t)$ pair, the problem is to maximize the number of observed $(s,t)$ pairs. 
It is also possible to take $u_l$ as different values for each $(s,t)$ pair.
For example, if $u_l$ is the OD traffic observed by the cut, the problem can be interpreted as maximizing the total observed traffic.
Thus, CSP2 can be formulated as a binary integer programming problem, which can also be solved with a standard solver.

\subsection{Properties on the upper bound of the solution of CSP1}\label{sec:char_cut}
The properties on the solution's upper bound of CSP1 are derived.
Let $\underline{d}(s)$ be the number of downstream links from the node with origin $s$.
This set of links is a cut since if all the downstream links from the origin are removed, then none of the endpoints can be reached.
This is called as the set of links $c_s^{st}$.
Similarly, let $\overline{d}(t)$ denote the number of all links inflowing to the destination $t$.
This set of links is also cut and is denoted by $c_t^{st}$.

First, the $(s,t)$-cut upper bound characteristics are shown in below.

\textbf{Lemma 1} \ The Eq.\eqref{eq:st_cut_lemma} holds for any $(s,t)$.
\begin{align}
  \label{eq:st_cut_lemma}
  \begin{split}
    \min_{c_m^{st}\in C^{st}}{\{ |c_m^{st}|\} } \leq \underline{d}(s) \quad \forall st|(s,t)\in W \\
    \min_{c_m^{st}\in C^{st}}{\{ |c_m^{st}|\} } \leq \overline{d}(t) \quad \forall st|(s,t)\in W
  \end{split}
\end{align}

\textbf{Proof.} \ By placing counters on the links that make up $c_s^{st}$, the traffic between $(s,t)$ is always observed.
Therefore, the $(s,t)$ minimum cut $\min_{c_m^{st}\in C^{st}}{\{ |c_m^{st}|\} }$ has an upper bound value of $\underline{d}(s)$.
The same is true for the links that make up $c_t^{st}$, and the Eq.\eqref{eq:st_cut_lemma} holds.

\rightline{Q.E.D.}

Then, other characteristics of the cut that separate the set of $(s,t)$ with the same origin or the same destination are introduced.
Let $f_{G(V,E)}(i,j)$ be a function indicating the existence of a directed path between $(i,j)$ in $G(V,E)$. If $f_{G(V,E)}(i,j)$ is one, a directed path exists in $(i,j)$; otherwise, it outputs zero.
For example, any $(s,t)$-cut $c_m^{st}$ does not have a directed path between $(s,t)$, and $f_{G(V,E-c_m^{st} )} (s,t)=0$ is guaranteed.
Let $C^{s*} = \bigcup_{t \in Q-s}{C^{st}}$ be the set of cuts that observe all $(s,t)$ with the same origin $s$.
Similarly, the set of cuts that observe all $(s,t)$ with the same destination $t$ is $C^{*t} = \bigcup_{s \in Q-t}{C^{st}}$.
There is an upper bound for the minimum number of links to observe all $(s,t)$ with a certain origin $s$.

\textbf{Lemma 2} \ The minimum number of links to observe between all $(s,t)$ with any starting point $s$ satisfies the Eq.\eqref{eq:s_cut_lemma}, and the number of links in a cut to observe between all $(s,t)$ with any destination $t$ satisfies the Eq.\eqref{eq:t_cut_lemma} is satisfied. 

\begin{align}
  \label{eq:s_cut_lemma}
  \begin{split}
    \min_{c_m^{st}\in C^{s*}}{\{|c_m^{st} \mid |  f_{G(V,E-c_m^{st})}(s,t)=0 \ \forall {t\in Q-s} \}} \leq \underline{d}(s) \\
  \end{split}
\end{align}
\begin{align}
  \label{eq:t_cut_lemma}
  \begin{split}
    \min_{c_m^{st}\in C^{*t}}{\{|c_m^{st} \mid |  f_{G(V,E-c_m^{st})}(s,t)=0 \ \forall {s\in Q-t} \}} \leq \overline{d}(t) \\
  \end{split}
\end{align}

\textbf{Proof.} \ By placing counters on the links that make up $c_s^{st}$, the traffic between all $(s,t)$ with a certain origin $s$ is always observed.
Therefore, the minimum cut $\min_{c_m^{st}\in C^{s*}}{\{ |c_m^{st}|\} }$ has an upper bound value of $\underline{d}(s)$.
The links that make $c_t^{st}$ are the same for the destination, and the Eq. \eqref{eq:s_cut_lemma} holds.

\rightline{Q.E.D.}

In road networks, the degree of the origin and destination nodes is often less than four, and the $(s,t)$-cut is often determined by the degree of one of the origin or destination nodes.
This is due to the diversity of the paths between the origin and destination nodes, which makes it difficult for cuts below the order of the origin and destination nodes to occur except for the links connecting to them.
\textbf{Lemma 2} shows that between single origin, multiple endpoints $(s,t)$, it is sufficient to place a sensor on the outflow link at the origin, or single destination, multiple origin $(s,t)$ at the end inflow link.
It is easy to imagine that in a road network, it would be difficult to find a combination of links that separates all multiple endpoints simultaneously with fewer than the number of outgoing links at the origin.
Therefore, \textbf{Lemma 2} gives a strong upper bound for the road network.
\textbf{Lemma 2} gave upper bound values for the number of links installed at each origin and at each destination.
Using this characteristic, the upper bound for the solution of CSP1 is also shown.

\textbf{Proposition 1} \ The optimal solution $Z^*$ of CSP1 satisfies the Eq.\eqref{eq:all_cut_prop}.
\begin{align}
  \label{eq:all_cut_prop}
  \begin{split}
    Z^* \leq \sum_{s\in Q}{\underline{d}(s)} \\
    Z^* \leq \sum_{t\in Q}{\overline{d}(t)} 
  \end{split}
\end{align}

\textbf{Proof.} \ If \textbf{Lemma 2} is applied to all origins and a counter is placed on $c_s^{st}$ at all origins $s$, then all $(s,t)$ between $(s,t)$ are observable.
Therefore, the optimal solution for CSP1 has an upper bound $\sum_{s\in Q}{\underline{d}(s)}$.
Similarly, if \textbf{Lemma 2} is applied to all destinations, the result is Eq.\eqref{eq:all_cut_prop}.

\rightline{Q.E.D.}

Although this property is an extension of the property shown in $\textbf{Lemma 2}$ to all origin and destination nodes, it is still expected to yield a solution equal to this upper bound, as shown in Chapter 4, since links connected to origin and destination nodes are often the most efficient to observe all $(s,t)$ pairs in a road network.
In addition, introducing the Eq. \eqref{eq:all_cut_prop} as a constraint in CSP1 can significantly limit the search region of $\mathbf{y}$, which is expected to improve the computational efficiency.

\subsection{A discussion of the cuts that should inputs into the CSP}

The method is shown in Section \textbf{\ref{sec:cut_enu}} enumerates all $(s,t)$-cuts. Depending on the network size, the amount of these cuts can be expected to be huge.
Introducing all $(s,t)$-cuts into CSP guarantees a global optimal solution. Therefore, if abundant computational resources are available, $\Psi$ can be generated by considering all $(s,t)$-cuts.
However, since the decision variable $\mathbf{y}$ for cut adoption is a binary integer in both CSP1 and CSP2, an increase in its size leads to an exponential increase in computational load.
However, for both CSP1 and CSP2, the objective function favors the selection of cuts with a small number of links, so it may be expected that optimal solutions can be obtained without considering cuts with too many links.
Note that, unlike the method of enumerating paths, the approach in this paper does not yield a solution that deviates from the OD covering rule, even if the variables input to the optimization problem are limited due to the nature of the cuts.
The solution obtained when the input cuts are limited is an approximate solution.
Based on the properties shown in Section \textbf{\ref{sec:char_cut}}, I will consider how to identify the cuts to be input to the CSP from the enumerated cuts.

In CSP1, since \textbf{Lemma 2} gives a strong upper bound, it would be a valid option to seek computational efficiency by restricting the input to CSP1 to only those cuts that satisfy Eq.\eqref{eq:s_cut_lemma} at each origin node.
However, the property shown in D is unavailable for CSP2, where the optimization problem specifies the $(s,t)$ to be observed.
At the time of this writing, I have not found a condition for the cut that should be entered in CSP2. Identifying these conditions is a subject for future work.
In this paper, I set an upper limit to the number of links that constitute a cut and input the corresponding cuts to CSP2.

\section{Numerical computation in the Sioux-Falls network}

\begin{figure}[t]
  \centering
  \includegraphics[keepaspectratio,scale=1.0]
  {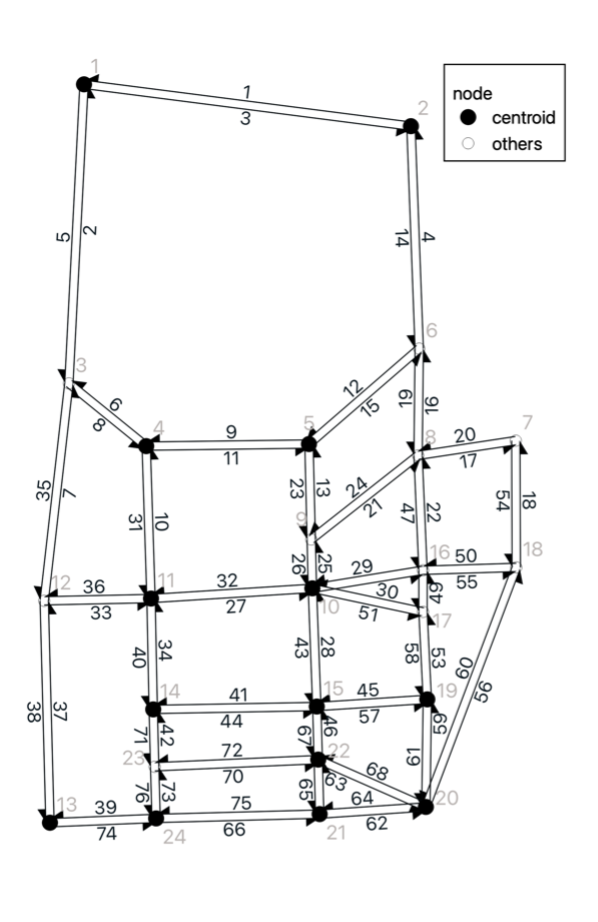}
  \caption{Sioux-Falls network}
  \label{fig:siouxfalls}
\end{figure}

This section demonstrates the effectiveness of the proposed method by applying it to the Sioux-Falls network.
The network conditions are the same as those in the application example of Yang et al.(2006), and the nodes shown in black in \figref{siouxfalls} are set as centroids.
The number of centroids is 14, and the number of $(s,t)$ pairs is 182 ($=14\times13$).

The numerical computation was performed on an Apple M2 Max 12 core CPU, 64 GB RAM, and Mac OS 13.4 OS and coded in Python.
\citet{provan_1996} was implemented using the \verb#st_all_cuts# function of the Python-igraph library.
Gurobi was used to solve CSP1 and CSP2.

\subsection{Enumerated cuts}

Enumeration of all $(s,t)$-cuts took less than 10 seconds. The resulting cuts are tabulated by the number of included links and shown in \tabref{cut_enumeration}.
After removing duplicates, the number of cuts for all link numbers is reduced, and the same links are included in many cuts for each $(s,t)$ pair.
The number of cuts enumerated for 12 links was the largest.

\begin{table}[t]
  \caption{Number of cuts enumerated in Sioux-Falls network}
  \label{tab:cut_enumeration}
  \centering
  \begin{tabular}{l l l}\hline
  Number of links & 
  \begin{tabular}{l}
    The number of cuts\\allowing OD duplication
  \end{tabular}&
  \begin{tabular}{l}
    The number of cuts\\without duplication
  \end {tabular} \\\hline
  2 & 126 & 8\\
  3 & 378 & 24\\
  4 & 1,088 & 52\\
  5 & 4,236 & 144\\ 
  6 & 12,976 & 370\\
  7 & 31,114 & 814\\
  8 & 66,168 & 1,656\\
  9 & 133,604 & 3,198\\
  10 & 254,234 & 5,838\\
  11 & 408,024 & 9,122\\
  12 & 508,776& 11,184\\
  13 & 491,842 & 10,662\\
  14 & 361,278 & 7,736\\
  15 & 217,320 & 4,602\\
  \hline
  \end{tabular}
\end{table}

\subsection{Optimal layout results in CSP1}
The cut consisting of the number of links that is less than or equal to the degree of all origins was used as the input for CSP1 according to the upper boundary values given in \textbf{Lemma 2}.
All $t_r$ were set to the same value to minimize the number of links to be installed.
The solution of CSP1 took less than one second, and 45 links were identified, which is consistent with the optimal solution in \citet{Yang2006}.
This optimal solution is consistent with the upper bound of the \textbf{Lemma 1}.
\begin{figure}[t]
  \centering
  \includegraphics[keepaspectratio,scale=1.0]
  {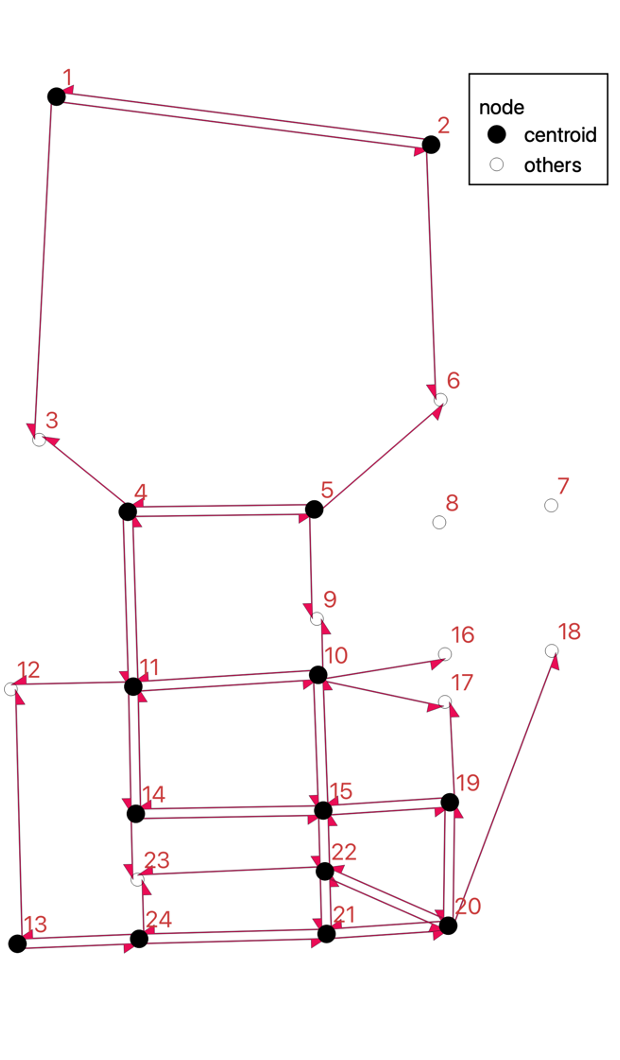}
  \caption{Layout of links with counters in the optimal solution}
  \label{fig:with_sensor}
\end{figure}
\begin{figure}[t]
  \centering
  \includegraphics[keepaspectratio,scale=1.0]
  {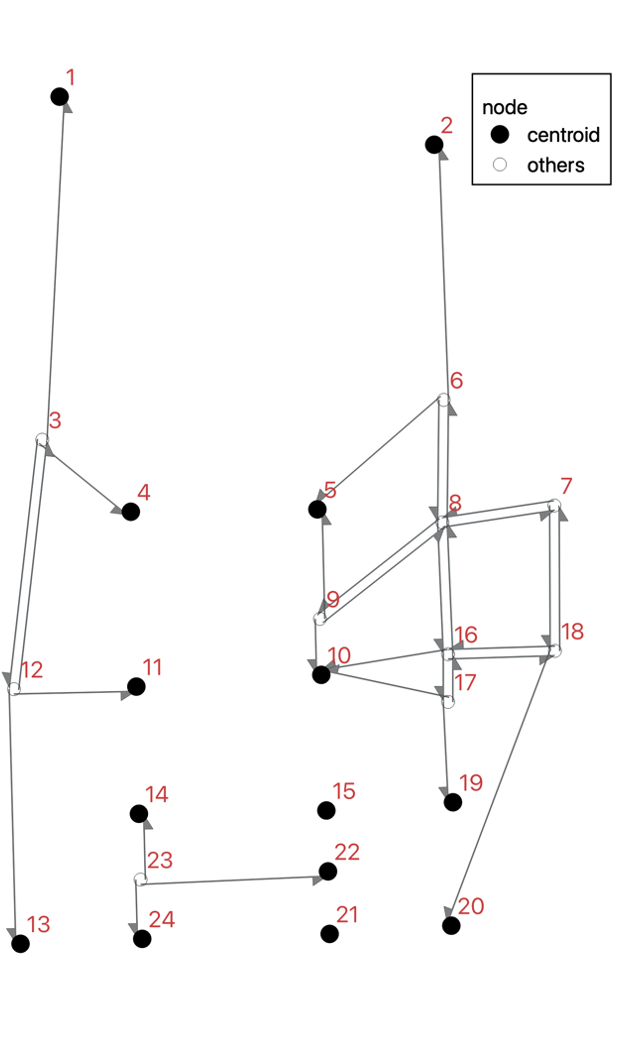}
  \caption{Distribution of links without counters}
  \label{fig:without_sensor}
\end{figure}

As can be seen in \figref{with_sensor}, all links originating from the centroids are selected in the optimal solution.
This means that the solution is consistent with the upper bound shown in \textbf{Lemma 2}.
Also, as can be confirmed by \figref{without_sensor}, there is no path from any of the centroids to the other centroids.
The solution is equivalent to the optimal solution presented in the previous literature and has been provided and confirmed to satisfy the OD Covering Rule.
\citet{Yang2006} do not give the computation time, but the row generation is repeated several hundred times.
On the other hand, our computation time of CSP1 is less than one second, and even when combined with cut enumeration, it takes only a dozen seconds.

\subsection{Optimal layout result in CSP2}

In CSP2, a budget for the sensor-located links must be provided.
As in \citet{Yang2006}, the cases with an upper bound on the number of links for every four links up to $K={4,8,12.... .48}$ as in the budget constraint are computed.
All enumerated cuts were input, and the calculation was attempted, but it did not end even after one hour.
Therefore, among the enumerated cuts, The cases where the cuts to be input to CSP2 were composed of a number of links according to multiple conditions based on the number of links included, such as {the degree of the origin, four or less, five or less, ... eight or less} are considered.
$u_l$ was the same for all cuts, and the problem was to maximize the number of observed $(s,t)$.

\begin{figure}[t]
  \centering
  \includegraphics[keepaspectratio,scale=0.4]
  {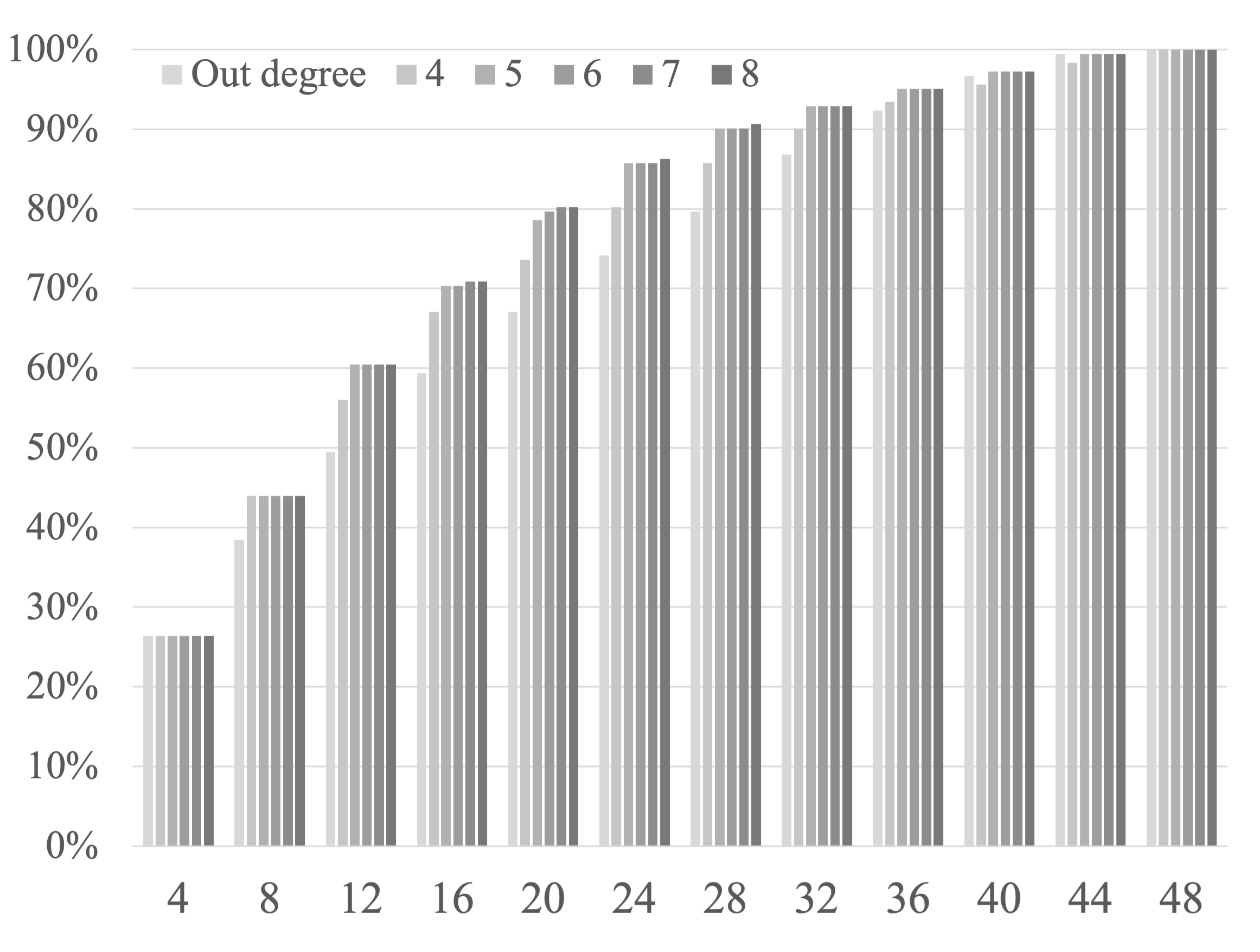}
  \caption{The ratio of the number of observed $(s,t)$ for each budget case and input cut condition}
  \label{fig:result_csp2}
\end{figure}

The number of observed $(s,t)$ for each case of the number of links to be placed and the conditions of the cuts to be input to CSP2 were calculated.
The obtained solutions were divided by 182, the total number of $(s,t)$, and the ratio of the number of $(s,t)$ that can be observed among all $(s,t)$ is shown in \figref{result_csp2}.

\subsubsection{Results on budget constraints}

The results were organized with a focus on the budget constraints.
Since the optimal solution of CSP1 is 45, all $(s,t)$ pairs can be observed if there are 45 links.
Therefore, when the budget constraint is set to 48 links, all $(s,t)$ pairs are observed.
When the budget constraint is 4, about 27\%  of the $(s,t)$ pairs are observed, but as the constraint is relaxed, the number of observed $(s,t)$ pairs increases.
However, when the budget constraint is 32 links, more than 90\% of the $(s,t)$ pairs are observed, and the increase in the number of observed $(s,t)$ pairs ratio due to the relaxation of the budget constraint is reduced.

\subsubsection{Summary of results regarding cuts and computation time to be input into CSP2}
Here, the results regarding the cuts to be input into CSP2 are summarized.
As mentioned above, when the budget constraint is small, the number of observed $(s,t)$ pairs is small, but as the budget constraint is relaxed, the number of observed $(s,t)$ pairs increases.
Two conditions are compared: one in which the input cuts are up to the same number of links as the degree of the starting point, and the other in which cuts consisting of up to eight links are input.
The results are the same for all budget cases when the budget constraint is four links, but the number of observed $(s,t)$ pairs is significantly different for budget constraints of 8 to 40.
There is no difference in the number of observed $(s,t)$ pairs for budget constraints of 40 or more.
This is because when the budget constraint is extremely small, according to the property of \textbf{Lemma 2}, selecting a cut consisting of a departing link of one starting point will result in the optimal solution because all $(s,t)$ pairs sharing that origin can be observed.
Also, when the budget constraint is large, the cut consisting of the departing link of the starting point is selected because the conditions of CSP1 and the problem are close.
In these conditions, even if only cuts consisting of links up to the degree are input, the results are such that the cuts consisting of the same number of links as the degree are adopted, so no difference is expected.
On the other hand, for budget cases between these conditions, there are cases where a cut can be shared by multiple $(s,t)$ pairs.
In those cases, a cut shared by multiple origins may be efficient even if it contains more links than the degree.
Examples of cuts that meet this condition will be discussed later.
From the above, it can be seen that it is better to input cuts consisting of a larger number of links, including those up to the degree unless a large number of $(s,t)$ pairs are observed compared to the total amount, or unless there is an extremely strong budget constraint.

\begin{figure}[t]
  \centering
  \includegraphics[keepaspectratio,scale=0.4]
  {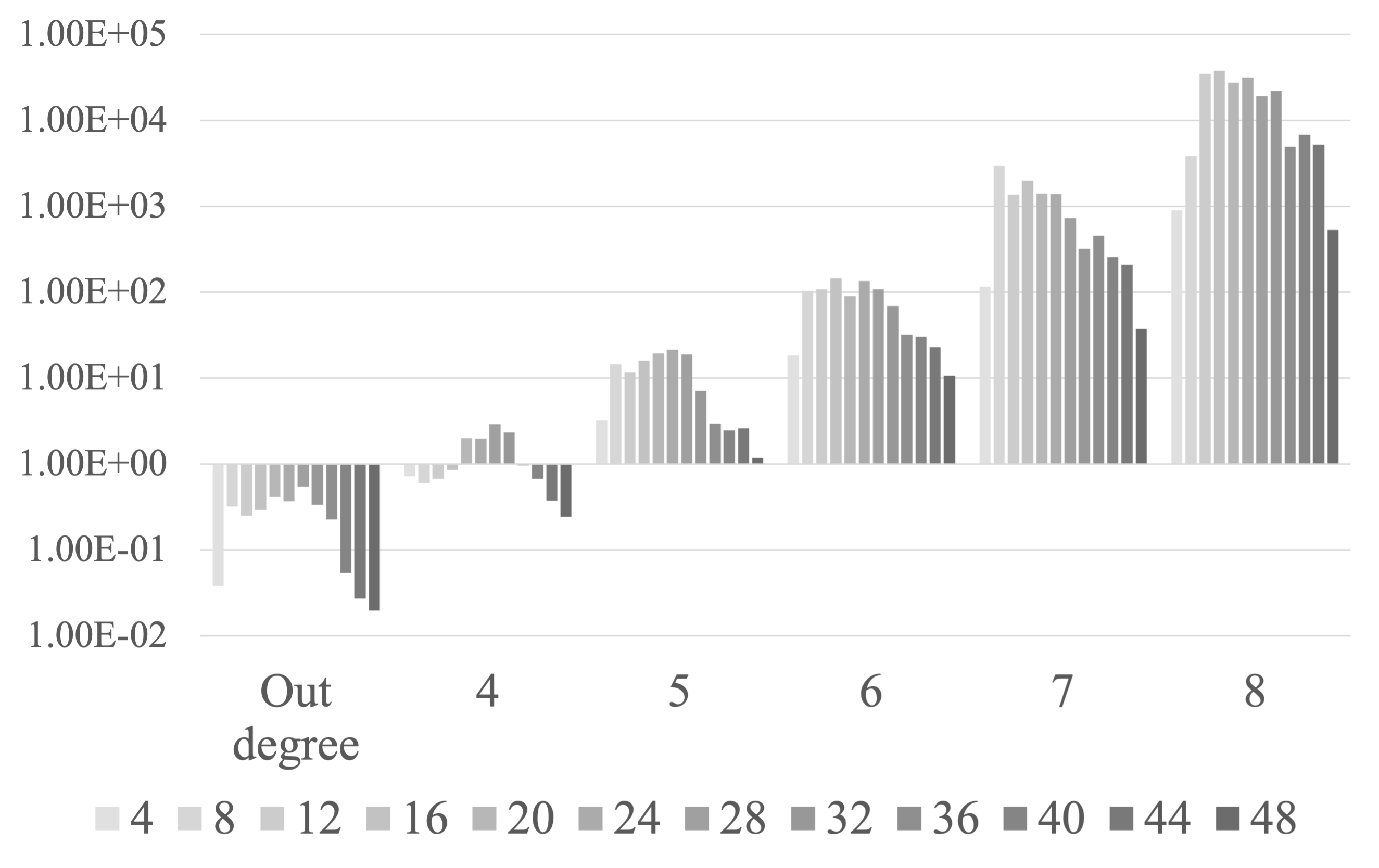}
  \caption{Observed $(s,t)$ pairs ratio for each budget case and input cut condition}
  \label{fig:comp_time_csp2}
\end{figure}

The computation time results for the cuts' conditions to be input to CSP2 are shown in \figref{comp_time_csp2}.

The horizontal axis shows each budget case, the bar's color shows the condition of the cuts to be input, and the vertical axis shows the time until the calculation is completed.
The vertical axis is a logarithmic scale, and the unit is seconds.
Focusing on the conditions of the cuts to be input, the calculation time is shorter when only the cuts with fewer links are considered for all budget cases.
Since the calculation time increases linearly on the logarithmic scale, it can be confirmed that the condition of the cuts to be input significantly impacts it.
If an approximate solution is acceptable for analytical purposes, it is better to consider only cuts with fewer links to obtain the solution quickly.

Furthermore, the computation time for each budget case was shorter for cases with large or small budget constraints and larger for budget cases in between.
The cause is unclear, but many solutions in these cases may be close to optimal.
In such cases, it may take some time for the upper and lower bound values to meet in Brunch and Bound.

\subsubsection{Cut shared among many $(s,t)$ pairs}

\begin{figure}[t]
  \centering
  \includegraphics[keepaspectratio,scale=0.55]
  {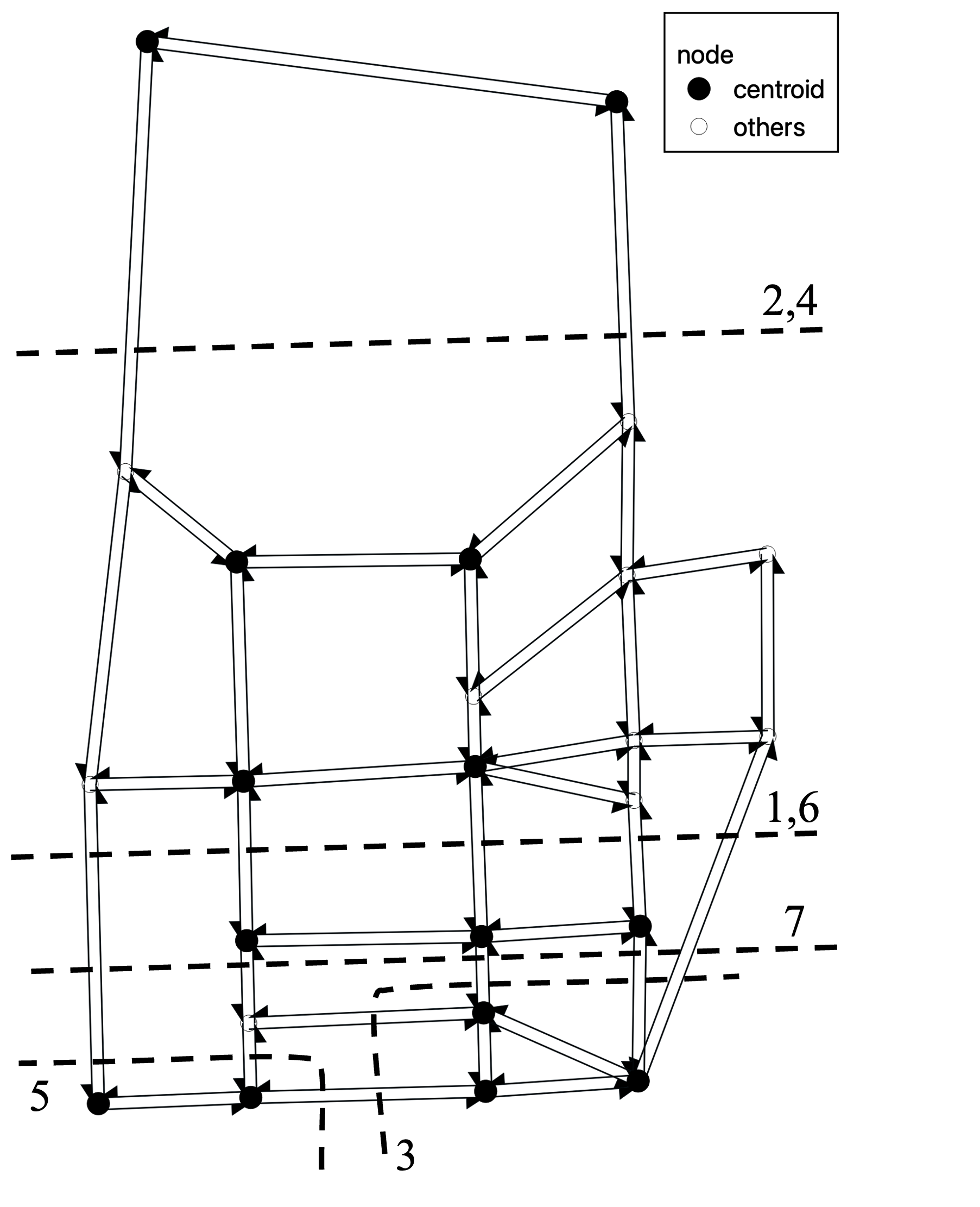}
  \caption{Layout of cuts shared among many $(s,t)$ pairs}
  \label{fig:shared_cuts}
\end{figure}

To confirm the presence or absence of cuts shared by many $(s,t)$ pairs, cuts with up to eight links are input into CSP2, and the budget constraint is set to $K=24$.
The cuts for which $y_l=1$ is the optimal solution are extracted, and the number of $(s,t)$ pairs that share the same link is counted.
The cuts that shared 10 or more $(s,t)$ pairs are shown in \tabref{shared_cuts}.
The corresponding layout is shown in \figref{shared_cuts}.
The dashed lines in the figure indicate the cuts, and the numbers above them are consistent with the Rank in \tabref{shared_cuts}.
If multiple numbers are assigned to the same cut, it indicates that the cut is composed of the same links but in different directions.

\begin{table}[t]
  \caption{Cut shared among many $(s,t)$}
  \label{tab:shared_cuts}
  \centering
  \begin{tabular}{l l l }\hline
  Rank & Link id 
    & \begin{tabular}{c} Number of\\
      shared $(s,t)$ pairs\end{tabular} \\\hline
  1 & 28,34,37,53,56 & 26\\
  2 & 2,4 & 24\\
  3 & 60,61,66,67,70 & 23\\
  4 & 5,14 & 20\\
  5 & 38,75,76 & 19\\
  6 & 38,40,43,58,60 & 18\\
  7 & 37,42,46,56,59 & 17\\
  \hline
  \end{tabular}
\end{table}

The results show, for example, that the cuts of Rank one, three, six, and seven each contain five links.
There is no centroid with a degree of five in the Sioux-Falls network.
Therefore, these cuts are not any of the $(s,t)$ minimum cuts but are selected as cuts that efficiently observe multiple ODs.
The number of shared $(s,t)$ pairs is 26 for Rank one cuts.
If the $c_s^{st}$ of a certain starting point is adopted, the number of $(s,t)$ that can be observed is 13.
Considering that $|c_s^{st}|$ is often two to four, it can be seen that the cuts of Rank one, three, six, and seven efficiently observe many $(s,t)$ pairs with five links.

\section{Conclusion}

This paper gives a formulation by cuts in SCLP, a counter location problem for efficiently observing traffic between $(s,t)$ with fewer links.
The proposed model showed that it can be formulated simply as an extension of the maximum weight closure problem.
Two problems corresponding to those given by \citet{Yang2006} were provided as CSP1 and CSP2, respectively.
This paper gains insight into the characteristics of the cut and shows that the sum of the degree of the origin or destination is a strong upper boundary value for CSP1.

The proposed method was verified for the Sioux-Falls network.
CSP1 was solved in less than one second in the author's environment, and the obtained solution was confirmed to match the optimal solution confirmed in previous studies.
Furthermore, it was confirmed that this solution matches the upper bound value of CSP1 shown in this paper.
In CSP2, it was confirmed that the number of $(s,t)$ pairs that can be observed decreases as the cuts to be input to CSP2 are limited according to the number of links.
Also, as the budget constraint is relaxed, the number of observed OD pairs increases, but as the upper limit of the number of links increases, the number of $(s,t)$ pairs that can be observed decreases.

This paper provides a formulation as an extension of MWCP, which can be solved faster than previous studies because it does not require repeated calculations of path search by column generation.
However, it requires enumerating all cuts between $(s,t)$ in advance.
In the Sioux-Falls network, it was possible to enumerate within 10 seconds in the author's environment, but calculation time may be a problem as the network becomes larger.
Also, since computation time increases as the number of cuts to be input to CSP2 increases, this issue must be addressed.
There are no studies where an exact solution method that satisfies the OD covering rule is applied in large networks in SCLP, but to ensure more general practicality, it is necessary to develop a methodology that can be solved more quickly.
One approach is to consider a selection method based on the optimal solution's properties for the size of the cuts to be input to CSP2.
Another approach is to construct a methodology that determines the cuts endogenously and determines the counter layout links in bulk as an optimization problem so that the cuts do not need to be enumerated in advance.

\section*{Acknowledgment}
This secondary publication is an English translation of the Japanese Journal of JSCE Transactions, Vol. 80, No. 5, 23-00244 (https://doi.org/10.2208/jscejj.23-00244).
This work was supported by JSPS KAKENHI Grant Numbers 21H01446 and 22H01610.

\bibliographystyle{unsrtnat}
\bibliography{doron_sugiura}  






\end{document}